\documentclass[11pt,a4paper]{article}

\usepackage{titlesec}
\usepackage{fancyhdr}
\usepackage{a4wide}
\usepackage{graphicx}
\usepackage{float}
\usepackage{amssymb}
\usepackage{amsmath}
\usepackage{amsthm}
\usepackage{color}
\usepackage{mathrsfs}
\usepackage{array}
\usepackage{eucal}
\usepackage{tikz}
\usepackage[T1]{fontenc}
\usepackage{inputenc}
\usepackage[english,french]{babel}
\usepackage{frcursive}
\usepackage{lmodern}
\usepackage{hyperref}
\usepackage{geometry}
\usepackage{changepage}
\changepage{0pt}{}{}{}{}{0pt}{}{0pt}{5pt}
\usepackage[numbers]{natbib}
\usepackage{hyperref}
\hypersetup{
pdfpagemode=none,
pdftoolbar=true,        
pdfmenubar=true,        
pdffitwindow=false,     
pdfstartview={Fit},    
pdftitle={Transmission eigenvalues with artificial background for explicit material index identification},    
pdfauthor={L. Audibert, L. Chesnel, H. Haddar},     
pdfsubject={},  
pdfcreator={L. Audibert, L. Chesnel, H. Haddar},   
pdfproducer={L. Audibert, L. Chesnel, H. Haddar}, 
pdfkeywords={}, 
pdfnewwindow=true,      
colorlinks=true,       
linkcolor=magenta,          
citecolor=red,        
filecolor=cyan,      
urlcolor=blue           
}

\newcommand{\dsp}{\displaystyle}

\newcommand{\Om}{\Omega}

\newcommand{\R}{\mathbb{R}}
\newcommand{\mL}{\mrm{L}}
\newcommand{\mH}{\mrm{H}}

\newcommand{\loc}{\mbox{\scriptsize loc}}
\newtheorem{theorem}{Theorem}[section]

\def\R{\mathbb R}

\def\Om{\Omega}
\def\Omb{\Omega_b}

\newcommand{\mrm}[1]{\mathrm{#1}}

\DeclareMathOperator{\supp}{\mathrm{supp}}
\newcommand{\Cplx}{\mathbb{C}}

\newcommand{\bfx}{\boldsymbol{x}}

\newcommand{\mZero}{\setminus \lbrace 0 \rbrace}
\newcommand{\Hinc}{\mrm{H}_{\mathrm{inc}}(\Om)}
\newcommand{\Hinct}{\tilde{\mrm{H}}_{\mathrm{inc}}(\Omb)}
\newcommand{\Hincd}{\tilde{\mrm{H}}_{\mathrm{inc}}(\partial \Omb)}

\begin{document}

~\vspace{0.4cm}
\begin{center}
{\sc \bf\LARGE 
\scalebox{1.04}{Transmission eigenvalues with artificial background}}\\[6pt]
{\sc \bf\LARGE 
\scalebox{1.04}{for explicit material index identification}}\\[15pt]
\end{center}

\begin{center}
\textsc{Lorenzo Audibert}$^{1,2}$, \textsc{Lucas Chesnel}$^2$, \textsc{Houssem Haddar}$^{2}$\\[18pt]
\begin{minipage}{0.9\textwidth}
{\small
$^1$ Department STEP, EDF R\&D, 6 quai Watier, 78401, Chatou CEDEX, France;\\[2pt]
$^2$ INRIA/Centre de mathématiques appliquées, \'Ecole Polytechnique, Université Paris-Saclay, Route de Saclay, 91128 Palaiseau, France.\\[10pt]
E-mails: \scalebox{0.9}{\texttt{lorenzo.audibert@edf.fr}, \texttt{lucas.chesnel@inria.fr}, \texttt{houssem.haddar@inria.fr}}\\[-14pt]
\begin{center}
(\today)
\end{center}
}
\end{minipage}
\end{center}
\vspace{0.4cm}

\noindent\textbf{Abstract.} 
We are interested in the problem of retrieving information on the refractive index $n$ of a penetrable inclusion embedded in a reference medium from farfield data associated with incident plane waves. Our approach relies on the use of  transmission eigenvalues (TEs) that carry information on $n$ and that can be determined from the knowledge of the farfield operator $F$. In this note, we explain how to modify $F$ into a farfield operator $F^{\mrm{art}}=F-\tilde{F}$, where $\tilde{F}$ is computed numerically, corresponding to well chosen artificial background and for which the associated TEs provide more accessible information on $n$.\\
\newline
\noindent\textbf{Keywords.}  Inverse scattering problems, linear sampling method, qualitative methods.\\
\newline
\noindent\textbf{Mathematics Subject Classification.}  35R30, 65M32, 35Q60, 35J40. \\

\selectlanguage{french}

\noindent\scalebox{1.05}{\textbf{Une identification d'indice explicite au moyen de valeurs propres de transmission}}\\
\noindent\scalebox{1.05}{\textbf{pour un milieu de référence artificiel}}\\[3pt]
\noindent\textbf{Résumé.} 
Nous souhaitons retrouver l'indice $n$ d'une inclusion pénétrable dans un milieu de référence connu à partir de la donnée de champs lointains associés à des ondes planes incidentes. Pour ce faire, nous utilisons les valeurs propres de transmission  (VPT) qui dépendent de $n$ et qui peuvent être déterminées à partir de l'opérateur de champ lointain $F$. Dans cette note, nous expliquons comment modifier $F$ en un opérateur de champ lointain  $F^{\mrm{art}}=F-\tilde{F}$, où $\tilde{F}$ est calculé numériquement, correspondant à un milieu de référence artificiel et pour lequel les VPT associées fournissent une information plus directe sur $n$. 
\noindent

\selectlanguage{english}

\section{Introduction}
In recent years Sampling Methods offered different perspectives in solving time harmonic inverse scattering problems \cite{CaCoHa16}. In addition to allow for a non iterative scheme to retrieve the support of inhomogeneities from multistatic data, these methods revealed the possibility to construct from the data a spectrum related to the material properties. This spectrum corresponds to the set of Transmission Eigenvalues (TEs) of the so-called  Interior Transmission Problem (ITP) (see (\ref{OriginalITEP})).  In the justification of sampling methods, substantial efforts have been made to prove discreteness of the set of TEs \cite{CaHa13} because most of these methods fail at frequencies corresponding to these values. However, since the work in \cite{CaCH10}, exploiting the failure of the reconstruction methods at TEs, it was proved that they can be determined from measured data and therefore can be exploited to infer information on the material properties. The determination of TEs from measured data has been improved using the framework of the Generalized Linear Sampling Method (GLSM) where exact knowledge of the support is no longer needed \cite{AuHa14, CaCoHa16}. See also \cite{KiLe13} for a different approach.

Under certain assumptions on $n$, the refractive index of the considered inhomogeneity appearing in Problem (\ref{PbChampTotalFreeSpaceIntrod}) below, it has been proved that there exist an infinite number of real positive TEs ($k^2>0$) \cite{CaGH10}. Note that in practice only real positive TEs are of interest because one can only play with real wavenumbers for measurements. This result of existence of real positive TEs is not obvious because the ITP (see equation after (\ref{OriginalITEP}) below) is quadratic in $k^2$ and it does not seem possible to see the spectrum of (\ref{OriginalITEP}) as the spectrum of a self-adjoint operator. In particular, in 1D situations, it has been established that complex TEs do exist. 
 
Although mathematically interesting, relying on transmission eigenvalues  to determine quantitative features on $n$ is difficult. The reason is twofold. First, information is lost in complex eigenvalues which cannot be measured in practice. Second it is difficult to establish sharp estimates for real TEs with respect to $n$ due to the complexity of the problem. In this note, we explain how to work with another farfield operator $F^{\mrm{art}}$ corresponding to an artificial background (reference medium) for which the associated TEs have a more direct connection with $n$. Put differently, working with $F^{\mrm{art}}$, our goal is to simplify the solution of the inverse spectral problem consisting in determining $n$ from the knowledge of real positive TEs. Important in the analysis is the fact that $F^{\mrm{art}}$ is given by the formula $F^{\mrm{art}}=F-\tilde{F}$ where $\tilde{F}$ can be obtained via a rather direct numerical computation. Therefore, in practice $F^{\mrm{art}}$ can also be considered as a data. Interestingly also, our approach does not require a priori knowledge of the exact support of the inhomogeneity. It is sufficient to know that the defect in the reference medium is located in a given bounded region.

Close to our study are the papers \cite{GiHa12,CCMM16, AuCaHa17}. In the first one, the authors reformulate the ITP as an eigenvalue problem for the material coefficient. In the second and third ones, it is explained how to identify $n$ from the knowledge of $F(k)-\tilde{F}(k,\gamma)$ at a single wavenumber $k$ and for a range of $\gamma$. Here $F(k)-\tilde{F}(k,\gamma)$ can be seen as the farfield operator corresponding to a background depending on an artificial parameter $\gamma$. In comparison with our approach, this method is interesting because it requires to know $F$ at a single wavenumber ($\tilde{F}(k,\gamma)$ can be computed numerically). However, the relation between associated TEs and $n$ is a bit more complex than in our case.

\section{Setting}

We assume that the propagation of waves in time harmonic regime in the reference medium $\R^d$, $d=2,3$, is governed by the Helmholtz equation $\Delta u+k^2 u=0$ with $k>0$ being the wavenumber. The localized perturbation in the reference medium is modeled by some bounded open set $\Om\subset\R^d$ with Lipschitz boundary $\partial\Om$ and a refractive index $n\in \mL^{\infty}(\R^d)$. We assume that $n$ is real valued, that $n=1$ in $\R^d\setminus\Om$ and that $\mrm{ess}\inf_{\Om}n$ is positive. The scattering of the incident plane wave $u_{i}(\cdot,\boldsymbol{\theta}_{i}):=e^{i k \boldsymbol{\theta}_{i}\cdot\boldsymbol{x}}$ of direction of propagation $\boldsymbol{\theta}_{i}\in\mathbb{S}^{d-1}$ by $\Om$ is described by the problem
\begin{equation}\label{PbChampTotalFreeSpaceIntrod}
\begin{array}{|rcll}
\multicolumn{4}{|l}{\mbox{Find }u=u_i+u_s\mbox{ such that}}\\[2pt]
\Delta u+k^2n\,u & = & 0 & \mbox{ in }\R^d,\\[5pt]
\multicolumn{4}{|c}{\dsp\lim_{r\to +\infty} r^{\frac{d-1}{2}}\left( \frac{\partial u_{s}}{\partial r}-ik u_{s} \right)  = 0,}
\end{array}
\end{equation}
with $u_{i}=u_{i}(\cdot,\boldsymbol{\theta}_{i})$. The last line of (\ref{PbChampTotalFreeSpaceIntrod}), where $r=|\boldsymbol{x}|$, is the Sommerfeld radiation condition and is assumed to hold uniformly with respect to $\boldsymbol{\theta}_{s} = \boldsymbol{x} / r$.  For all $k>0$, Problem (\ref{PbChampTotalFreeSpaceIntrod}) has a unique solution  $u \in \mH^2_{\loc}(\R^d)$. The scattered field $u_{s}(\cdot,\boldsymbol{\theta}_{i})$ has the expansion 
\begin{equation}\label{scatteredFieldFreeSpaceIntro}
u_{s}(\boldsymbol{x},\boldsymbol{\theta}_{i})= \dsp e^{i k r}r^{-\frac{d-1}{2}}\,\Big(\,u_{s}^{\infty}(\boldsymbol{\theta}_{s},\boldsymbol{\theta}_{i})+O(1/r)\,\Big),
\end{equation}
as $r\to+\infty$, uniformly in $\boldsymbol{\theta}_{s}\in \mathbb{S}^{d-1}$.  The function $u^{\infty}_{s}(\cdot,\boldsymbol{\theta}_{i}): \mathbb{S}^{d-1}\to\Cplx$,  is called the farfield pattern associated with $u_{i}(\cdot,\boldsymbol{\theta}_i)$. From the farfield pattern, we can define the farfield operator $F:\mL^2(\mathbb{S}^{d-1})\to\mL^2(\mathbb{S}^{d-1})$ such that 
\begin{equation}\label{defInitialFFoperator}
(Fg)(\boldsymbol{\theta}_{s})=\int_{\mathbb{S}^{d-1}}g(\boldsymbol{\theta}_{i})\,u_{s}^{\infty}(\boldsymbol{\theta}_{s},\boldsymbol{\theta}_{i})\,ds(\boldsymbol{\theta}_{i}).
\end{equation}
The function $Fg$ corresponds to the farfield pattern for the scattered field in (\ref{PbChampTotalFreeSpaceIntrod}) with $u_i=u_i(g):=\int_{\mathbb{S}^{d-1}}g(\boldsymbol{\theta}_{i})e^{ik\boldsymbol{\theta}_{i}\cdot\bfx}\,ds(\boldsymbol{\theta}_{i})$ (Herglotz wave function). Define the operator $\mathcal{H}: \mL^2(\mathbb{S}^{d-1}) \to \mL^2(\Om)$ such that $\mathcal{H}g = u_i(g)|_\Om$ and the space $\Hinc :=\{v \in \mL^2(\Om); \; \Delta v+k^2 v = 0 \mbox{ in } \Om\}$. It is known that $\Hinc $ is nothing but the closure of the range of the operator $\mathcal{H}$ in $\mL^2(\Om)$. Observing that $\Delta u_s
+k^2nu_s=k^2(1-n)u_i(g)$ (in particular $u_s$ depends only on the values of $u_i(g)|_{\Om}$), we can factorize $F$ as $F=G\mathcal{H}$ where the operator $G:\Hinc\to\mL^2(\mathbb{S}^{d-1})$ is the extension by continuity of the mapping $u_i(g)|_{\Om}\mapsto u_s^{\infty}$. The real transmission eigenvalues are defined as the values of $k\in \R$ for which $G$ is not injective. In such a case, there is a (generalized) incident wave $v \in \Hinc$ such that the associated farfield is zero and therefore, by the Rellich Lemma, the scattered field $u_s$ is zero outside $\Om$. This leads to the equivalent definition of TEs as the values of $k\in \R$ for which the problem 
\begin{equation}\label{OriginalITEP}
\begin{array}{|rcll}
\Delta u_s+k^2 nu_s&=& k^2(1-n) v&\mbox{ in }\Om \\
\Delta v+k^2v&=&0&\mbox{ in }\Om \\
\end{array}
\end{equation}
admits a non trivial solution $(u_s,v)\in\mH^2_0(\Om) \times\mL^2(\Om)$. In particular, if $(n-1)^{-1}\in\mL^{\infty}(\Om)$, (\ref{OriginalITEP}) can be equivalently written as $(\Delta+k^2)((n-1)^{-1}(\Delta u_s+k^2 nu_s))=0 $ in $\Om$ with $u_s\in\mH^2_0(\Om)$. A method has been designed in \cite{l-thesis,CaCoHa16} to identify TEs from the farfield operator $F$, that we shall generalize later for a modified background. One of the main troubles with TEs is that their link with the index of refraction $n$ is not explicit nor easily accessible. Some monotonicity results have been obtained in \cite{CaGH10} but only for some of the TEs. We hereafter explain how a simple modification of the farfield operator leads to simpler transmission eigenvalue problems and more accessible information on the index of refraction. This idea was motivated by recent works on so-called Steklov eigenvalues \cite{CCMM16} and modified backgrounds with metamaterials \cite{AuCaHa17}.

\section{Transmission eigenvalues with ZIM background}
\label{sectionZim}
Assume that one has a priori knowledge of a Lipschitz domain $\Omb$ such that $\Om\subset\Omb$. We emphasize that we do not require to know exactly the support $\overline{\Om}$ of the defect in the reference medium. Consider the scattering problem 
\begin{equation}\label{PbChampTotalFreeSpaceComp}
\begin{array}{|rcll}
\multicolumn{4}{|l}{\mbox{Find }\tilde{u}=u_i+\tilde{u}_s\mbox{ such that}}\\[2pt]
\Delta\tilde{u}+k^2\rho\,\tilde{u} & = & 0 & \mbox{ in }\R^d,\\[5pt]
\multicolumn{4}{|c}{\dsp\lim_{r\to +\infty} r^{\frac{d-1}{2}}\left( \frac{\partial \tilde{u}_{s}}{\partial r}-ik \tilde{u}_{s} \right)  = 0}
\end{array}
\end{equation}
where $\rho$ is the function such that $\rho=0$ in $\Omb$ and $\rho=1$ in $\R^d\setminus\overline{\Omb}$. This artificial media is referred to as Zero-Index Material because we choose $\rho=0$ inside $\Omb$. This choice greatly simplifies the structure of the associated interior transmission problem as we shall see. For $u_i=u_i(\cdot,\boldsymbol{\theta}_{i})=e^{i k \boldsymbol{\theta}_{i}\cdot\boldsymbol{x}}$ with $\boldsymbol{\theta}_{i}\in\mathbb{S}^{d-1}$, as for (\ref{PbChampTotalFreeSpaceIntrod}), this problem admits a unique solution in $\mH^1_{\loc}(\R^d)$. We denote $\tilde{u}_{s}(\cdot,\boldsymbol{\theta}_{i})$ the associated scattered field and $\tilde{u}^{\infty}_{s}(\cdot,\boldsymbol{\theta}_{i}): \mathbb{S}^{d-1}\to\Cplx$ the corresponding farfield pattern. From the farfield pattern, we define the farfield operator $\tilde{F}:\mL^2(\mathbb{S}^{d-1})\to\mL^2(\mathbb{S}^{d-1})$ such that $(\tilde{F}g)(\boldsymbol{\theta}_{s})=\int_{\mathbb{S}^{d-1}}g(\boldsymbol{\theta}_{i})\,\tilde{u}_{s}^{\infty}(\boldsymbol{\theta}_{s},\boldsymbol{\theta}_{i})\,ds(\boldsymbol{\theta}_{i})$. Finally we define the artificial farfield operator $F^{\mrm{art}}:\mL^2(\mathbb{S}^{d-1})\to\mL^2(\mathbb{S}^{d-1})$ as
\begin{equation}\label{DefArtFFop}
F^{\mrm{art}}:=F-\tilde{F}.
\end{equation}
From a practical point of view, notice that $F$ is given by the measurements while $\tilde{F}$ has to be computed, which is achievable because Problem (\ref{PbChampTotalFreeSpaceComp}) does not involve $n$ (and $\Omb$ is known). Therefore, we can consider $F^{\mrm{art}}$ as a data for the inverse problem of determining $n$. Let us denote
\begin{equation}\label{DefuiArt}
u_i^{\mrm{art}}(g):=\int_{\mathbb{S}^{d-1}}g(\boldsymbol{\theta}_{i})\tilde{u}(,\boldsymbol{\theta}_{i}) \,ds(\boldsymbol{\theta}_{i})
\end{equation}
for given $g \in \mL^2({\mathbb{S}^{d-1}})$ and define the operator $\mathcal{H}^{\mrm{art}}: \mL^2(\mathbb{S}^{d-1}) \to \mL^2(\Omb)$ such that  $\mathcal{H}^{\mrm{art}}g = u_i^{\mrm{art}}(g)|_{\Omb}$. Set also $\Hinct :=\{v \in\mL^2(\Omb); \; \Delta v = 0 \mbox{ in } \Omb\}$. Observing that $\Delta(u_s-\tilde{u}_s)+k^2n(u_s-\tilde{u}_s)=k^2(\rho-n)u_i^{\mrm{art}}(g)$, we can factorize $F^{\mrm{art}}$ as $F^{\mrm{art}}=G^{\mrm{art}}\mathcal{H}^{\mrm{art}}$ where the operator $G^{\mrm{art}}:\Hinct\to\mL^2(\mathbb{S}^{d-1})$ is the extension by continuity of the mapping $u_i^{\mrm{art}}(g)|_{\Omb}\mapsto u_s^{\infty}-\tilde u_s^{\infty}$. Similarly as above, we now define the transmission eigenvalues as the values of $k$ for which the operator $G^{\mrm{art}}$ is not injective. Denoting $w:=u_s-\tilde u_s$, this is now equivalent, due to the Rellich Lemma, to define the transmission eigenvalues as the values of $k$ for which the problem 
\begin{equation}\label{eigenpb}
\begin{array}{|rcll}
\Delta w +k^2 nw&=& -k^2 n v&\mbox{ in }\Omb \\
\Delta v&=&0&\mbox{ in }\Omb 
\end{array}
\end{equation}
admits a non trivial solution $(w,v)\in \mH^2_0(\Omb) \times  \mL^2(\Omb)$. We observe that $k=0$ is an eigenvalue of infinite multiplicity of (\ref{eigenpb}). Now we consider the case $k\ne0$. Then we find that (\ref{eigenpb}) admits a non trivial solution if and only if there is $w\not\equiv0$ such that $(\Delta +k^2 n)(n^{-1}\Delta w)=0$ in $\Omb$, that is if and only if $(w,k)$ is a solution of the problem
\begin{equation}\label{eigenpbSndOrder}
\begin{array}{|l}
\mbox{Find } (w,k) \in \mH^2_0(\Omb) \mZero\times \R \mbox{  such that }\\[3pt]
\Delta(n^{-1}\Delta w)=-k^2\Delta w\mbox{ in }\Omb.
\end{array}
\end{equation}
In opposition with problem (\ref{OriginalITEP}) (more precisely, see equation after (\ref{OriginalITEP})), one ends up here with a well known linear eigenvalue problem similar to the so-called plate buckling eigenvalue problem. Classical results concerning linear self-adjoint compact operators guarantee that the spectrum of (\ref{eigenpbSndOrder}) is made of real positive isolated eigenvalues of finite multiplicity $0<\lambda_0 \le \lambda_1 \le \dots \le \lambda_p \le\dots$ (the numbering is chosen so that each eigenvalue is repeated according to its multiplicity). Moreover, there holds $\lim_{p\to+\infty}\lambda_p=+\infty$ and we have the \textit{min-max} formulas
\begin{equation}\label{FormuleMinMax}
\lambda_p=\min_{E_p\in\mathscr{E}_p}\max_{w\in E_p\setminus\{0\}}\cfrac{(n^{-1}\Delta w,\Delta w)_{\mL^2(\Omb)}}{\|\nabla w\|^2_{\mL^2(\Omb)}}.
\end{equation}
Here $\mathscr{E}_p$ denotes the sets of subspaces $E_p$ of $\mH^2_0(\Omb)$ of dimension $p$. Observe that the characterisation of the spectrum of Problem (\ref{eigenpb}) is much simpler than the one of Problem (\ref{OriginalITEP}). Moreover, it holds under very general assumptions for $n$: we just require that $n|_{\Omb}\in\mL^{\infty}(\Omb)$ with $\mrm{ess}\inf_{\Om}n>0$.  In particular, $n$ can be equal to one inside $\Omb$ (thus, as already mentioned, we do not need to know exactly the support $\overline{\Om}$ of the defect in the reference medium) and $n-1$ can change sign on the boundary. In comparison, the analysis of the spectrum of (\ref{OriginalITEP}) in such situations is much more complex and the functional framework must be adapted according to the values of $n$ (see \cite{Ches16} in the case where $n-1$ changes sign on $\partial\Om$). The second advantage of considering Problem (\ref{eigenpb}) instead of Problem (\ref{OriginalITEP}) is that the spectrum of (\ref{eigenpb}) is entirely real. This result is interesting combined with the following theorem. Let us set $\tilde{F}_\sharp := |\tilde{F} + \tilde{F}^{\ast}| +  |\tilde{F} - \tilde{F}^{\ast}| $, where $\tilde{F}^{\ast}$ is the adjoint of $\tilde{F}$, and define for $\alpha >0 $, $g$, $\phi\in\mL^2(\mathbb{S}^{d-1})$, the functional
\[
J_\alpha(g, \phi) := \alpha (\tilde F_\sharp g, g)_{\mL^2(\mathbb{S}^{d-1}) } + \| F^{\mrm{art}} g - \phi\|^2_{\mL^2(\mathbb{S}^{d-1}) }.
\]
Notice that in the penalty term, we use the operator $\tilde F_\sharp$ and not $F^{\mrm{art}}_\sharp$ defined similarly as $\tilde F_\sharp$. 
The reason is that we do not know if $g\mapsto(F^{\mrm{art}}_\sharp g, g)_{\mL^2(\mathbb{S}^{d-1}) }$ is equivalent to $\|u_i(g)\|_{\mL^2(\Omb)}^2$. We then consider for ${\boldsymbol{z}}\in \R^d$ a function $g_{{\boldsymbol{z}}}^\alpha\in \mL^2(\mathbb{S}^{d-1}) $ such that
\[
J_\alpha(g_{{\boldsymbol{z}}}^\alpha, \phi_{\boldsymbol{z}}^\infty) \le \alpha + \inf_{g \in \mL^2(\mathbb{S}^{d-1})} J_\alpha(g, \phi_{\boldsymbol{z}}^\infty)
\]
where $\phi_{\boldsymbol{z}}^\infty(\boldsymbol{\theta}_{s}) := e^{i k \boldsymbol{\theta}_{s}\cdot{\boldsymbol{z}}}$ is the farfield associated with a point source at the point ${\boldsymbol{z}}$.
\begin{theorem}\label{ThmCaractSpectrum}
Assume that  the farfield operator $F^{\mrm{art}}$ has dense range. Then $k^2$ is an eigenvalue of (\ref{eigenpbSndOrder}) if and only if the set of points ${\boldsymbol{z}}$ for which $ (\tilde F_\sharp g_{{\boldsymbol{z}}}^\alpha, g_{{\boldsymbol{z}}}^\alpha)_{\mL^2(\mathbb{S}^{d-1}) }$ is bounded as $\alpha \to 0$ is nowhere dense in $\Omb$. 
\end{theorem}
\noindent This theorem guarantees that peaks in the curve 
\begin{equation}\label{indicatorGLSM}
k\mapsto I(k,n) := \int_{\Omega_0 \subset \Omb} (\tilde F_\sharp g_{{\boldsymbol{z}}}^\alpha, g_{{\boldsymbol{z}}}^\alpha)_{\mL^2(\mathbb{S}^{d-1}) } \, d{\boldsymbol{z}}
\end{equation}
for small values of $\alpha$ correspond to $k^2$ which are eigenvalues of (\ref{eigenpb}). Therefore, from the knowledge of $F^{\mrm{art}}$ for real $k\in(0;+\infty)$, we can identify all the spectrum of (\ref{eigenpb}). In other words, there is no loss of information in complex eigenvalues which can occur for Problem (\ref{OriginalITEP}). The proof of Theorem \ref{ThmCaractSpectrum} is similar to the one of \cite[Theorem 7]{AuCaHa17} and uses the fact that $(\tilde F_\sharp g, g)_{\mL^2(\mathbb{S}^{d-1}) }$ is equivalent  to $\|u_i(g)\|_{\mL^2(\Omb)}^2$.  The latter can be deduced for instance from \cite[Lemma 2.33 and Theorem 2.31]{CaCoHa16}. We remark that the farfield operator $F^{\mrm{art}}$ fails to be of dense range only if $k^2$ is an eigenvalue of (\ref{eigenpbSndOrder}) such that the corresponding eigenfunction  $v$ in (\ref{eigenpb}) is of the form $v =u_i(g)=\mathcal{H}^{\mrm{art}}g$ for some $g$ in $\mL^2(\mathbb{S}^{d-1})$. Since $\mathcal{H}^{\mrm{art}}$ is a compact operator, cases where $v = \mathcal{H}^{\mrm{art}}g$ are in general exceptional. It has been shown in \cite{AlPaSy14}, in the case of the homogeneous background ($\rho\equiv1$ in $\R^d$), that this never happens for certain scattering objects with corners.\\
\newline 
Formula (\ref{FormuleMinMax}) is interesting because it guarantees that the $\lambda_p$ have monotonous dependence with respect to $n$. If we denote $\hat{\lambda}_p(\hat{n})$ the $\textit{p}-$th eigenvalue of Problem (\ref{eigenpbSndOrder}) with $n$ replaced by $\hat{n}$, then we have $\lambda_p(n)\le\hat{\lambda}_p(\hat{n})$ if $n\ge\hat{n}$. This estimate can be useful to derive qualitative information from a reference situation: for a given setting where the $\lambda_p(n)$ are known, if $n$ is changed into $\hat{n}$, one can have an idea of the nature of the perturbation.  From (\ref{FormuleMinMax}), we can also write
\[
\mrm{ess}\inf_{\Omb}n\le \lambda_p(1)/\lambda_p(n)\qquad\mbox{ and }\qquad \mrm{ess}\sup_{\Omb}n\ge \lambda_p(1)/\lambda_p(n).
\]
Note that the $\lambda_p(1)$ can be explicitly computed with a numerical code because they do not depend on $n$ and $\Omb$ is known. In particular, when $n$ is constant and $\Omb=\Om$, we have the formula
\begin{equation}\label{totoratio}
n=\lambda_p(1)/\lambda_p(n).
\end{equation}
Thus, in this case, from the knowledge of the farfield operator $F^{\mrm{art}}$ on an interval of wavenumbers $k$ containing an eigenvalue of Problem (\ref{eigenpbSndOrder}), we can identity the value of $n$. 

\section{A modified background with a non penetrable obstacle} \label{jevaismanger}
We discuss here other possibilities to determine simple spectral signatures associated with the index of refraction by modifying the background. More specifically we consider a background with a non penetrable obstacle $\Omb$ such that $\Om\subset\Omb$. On $\partial\Omb$, we prescribe the boundary condition $B(\tilde{u})=0$ where $B$ is a given boundary operator.  The construction is similar to above where (\ref{PbChampTotalFreeSpaceComp}) is replaced by 
\begin{equation}\label{PbChampTotalFreeSpaceComp2}
\begin{array}{|rlcl}
\Delta\tilde{u}+k^2\tilde{u}  &=&  0 &\mbox{ in }\R^d \setminus \overline\Omb
\\[5pt]
B(\tilde{u})  & = &0& \mbox{ on }\partial\Omb.
\end{array}
\end{equation}
We investigate here only the cases $B(\tilde{u}) =\tilde{u}$ and $B(\tilde{u}) = \partial \tilde{u} / \partial \nu + \gamma \tilde{u} $ which respectively correspond to Dirichlet and Robin scattering problems, with $\gamma \in \R$ being a fixed impedance parameter. The normal unit vector $\nu$ to $\partial\Omb$ is chosen directed to the exterior of $\Omb$. We redefine respectively $F^{\mrm{art}}$, $u_i^{\mrm{art}}(g)$ as in (\ref{DefArtFFop}), (\ref{DefuiArt}) from $\tilde{u}$ introduced in (\ref{PbChampTotalFreeSpaceComp2}). Then we consider the operator $\mathcal{H}^{\mrm{art}}: \mL^2(\mathbb{S}^{d-1}) \to \Hincd$ such that $\mathcal{H}^{\mrm{art}}g = B^*(u_i^{\mrm{art}}(g))|_{\partial \Omb}$ where $B^*(\varphi) =\partial\varphi / \partial \nu$ for  the Dirichlet problem and $B^*(\varphi) =\varphi - \gamma \partial \varphi / \partial \nu$ for the Robin problem. Here we take $\Hincd=\mH^{1/2}(\partial \Omb)$ for $\gamma =0$ and $\Hincd=\mH^{-1/2}(\partial \Omb)$ for other cases. Introduce the function $w \in \mH^1_{\loc}(\R^d\setminus\overline{\Omb}) \times \mH^1(\Omb) $ such that $w=u_s-\tilde{u}_s$ in $\R^d\setminus\overline{\Omb}$ and $w=u$ in $\Omb$. Note that $w$ satisfies the problem 
\begin{equation}\label{PbChampTotalFreeSpaceComp3}
\begin{array}{|rcll}
\Delta w + k^2 w &=& 0&\mbox{ in }\R^d\setminus \Omb\\[2pt]
\Delta w+k^2n\,w &= & 0 & \mbox{ in }\Omb\\[2pt]
\multicolumn{4}{|l}{\left[B(w)\right] = 0 \mbox{ and }\left[B^*(w)\right] =-\psi \mbox{ on }\partial\Omb}
\\[5pt]
\multicolumn{4}{|l}{\dsp\lim_{r\to +\infty} r^{\frac{d-1}{2}}\left( \frac{\partial w}{\partial r}-ik w\right)  = 0,}
\end{array}
\end{equation}
with $\psi= B^*(u_i^{\mrm{art}}(g))$ (here $[\cdot]$ denotes the jump across $\partial\Omb$ with respect to the orientation of the normal $\nu$). From this observation, we see that we can factorize $F^{\mrm{art}}$ as $F^{\mrm{art}}=G^{\mrm{art}}\mathcal{H}^{\mrm{art}}$ where the operator $G^{\mrm{art}}:\Hincd\to\mL^2(\mathbb{S}^{d-1})$ is the mapping $\psi\mapsto w^\infty$, $w^\infty$ being the farfield associated with the solution $w$ of problem (\ref{PbChampTotalFreeSpaceComp3}). Defining the transmission eigenvalues as the values of $k$ for which the operator $G^{\mrm{art}}$ is not injective, one ends up equivalently with $k^2$ being the eigenvalues of the cavity problem, $w \in\mH^1(\Omb)$
\begin{equation}\label{eigenpbobs}
\begin{array}{|rcll}
\Delta w +k^2 nw&=& 0&\mbox{ in }\Omb \\
B(w)&=&0&\mbox{ on }\partial \Omb.  
\end{array}
\end{equation}
One can prove the following identification theorem.
\begin{theorem}\label{ThmCaractSpectrum1}
Assume that  the farfield operator $F^{\mrm{art}}$ has dense range and that  $k^2$ is not an eigenvalue of (\ref{eigenpbobs}) for $n=1$. Then  $k^2$ is  an eigenvalue of  (\ref{eigenpbobs}) if and only if the set of points ${\boldsymbol{z}}$ for which $ (\tilde F_\sharp g_{{\boldsymbol{z}}}^\alpha, g_{{\boldsymbol{z}}}^\alpha)_{\mL^2(\mathbb{S}^{d-1}) }$ is bounded as $\alpha \to 0$ is nowhere dense in $\Omb$. 
\end{theorem}
\noindent The assumption on $k^2$ is required so that  $(\tilde F_\sharp g, g)_{\mL^2(\mathbb{S}^{d-1}) }$ is  equivalent  to $\|u_i(g)\|_{\Hincd}^2$ (see for instance \cite{KiGr08}). The proof is analogous to the one of Theorem \ref{ThmCaractSpectrum}. Then we obtain similar conclusions as in the previous section on the determination of $n$ from TEs. One has also the possibility to use all these simple spectra at once to determine $n$.

\section{Discussion of other choices of artificial backgrounds}
Note that for the artificial background in (\ref{PbChampTotalFreeSpaceComp}), we are not obliged to choose $\rho$ such that $\rho=0$ in $\Omb$ and $\rho=1$ in $\R^d\setminus\overline{\Omb}$. Consider an arbitrary real valued $\rho\in\mL^{\infty}(\R^d)$ such that $\rho-1$ is compactly supported. Let $\Omb$ be a domain such that $\supp(n-\rho)\subset\overline{\Omb}$. For such a $\rho$, setting $w:=u-v$, the corresponding transmission eigenvalues reads 
\[
\begin{array}{|rcll}
\Delta w+k^2 nw&=&k^2(\rho-n)v&\mbox{ in }\Omb \\
\Delta v+k^2 \rho v&=&0&\mbox{ in }\Omb 
\end{array}
\]
with $(w,v)\in\mH^2_0(\Om_b)\times\mL^2(\Om_b)$. And if $\rho$ is such that $n-\rho\ne0$ in $\Omb$, this leads to the problem
\begin{equation}\label{eigenpb3}
\begin{array}{|l}
\mbox{Find } (w,k^2) \in \mH^2_0(\Omb) \mZero\times \mathbb{C} \mbox{  such that }\\[3pt]
(\Delta+k^2\rho)\Big(\cfrac{1}{n-\rho}\,(\Delta w+k^2n w)\Big)=0\mbox{ in }\Omb.
\end{array}
\end{equation}
All the game with these artificial backgrounds consists in choosing $\rho$ such that TEs of (\ref{eigenpb3}) give interesting or simply usable information on $n$. 
One can also vary the boundary condition in the case discussed in Section \ref{jevaismanger} by varying the Robin parameter $\gamma$.  The natural question then would be to ask whether these spectra provide unique determination of the refractive index. This will be discussed in a future work.
 
\section{Numerical illustrations}\label{numSpectrum}
\begin{figure}[!ht]
\centering
\begin{tabular}{c}
\includegraphics[width=0.97\textwidth,trim={0 1.8cm 0 1.8cm}]{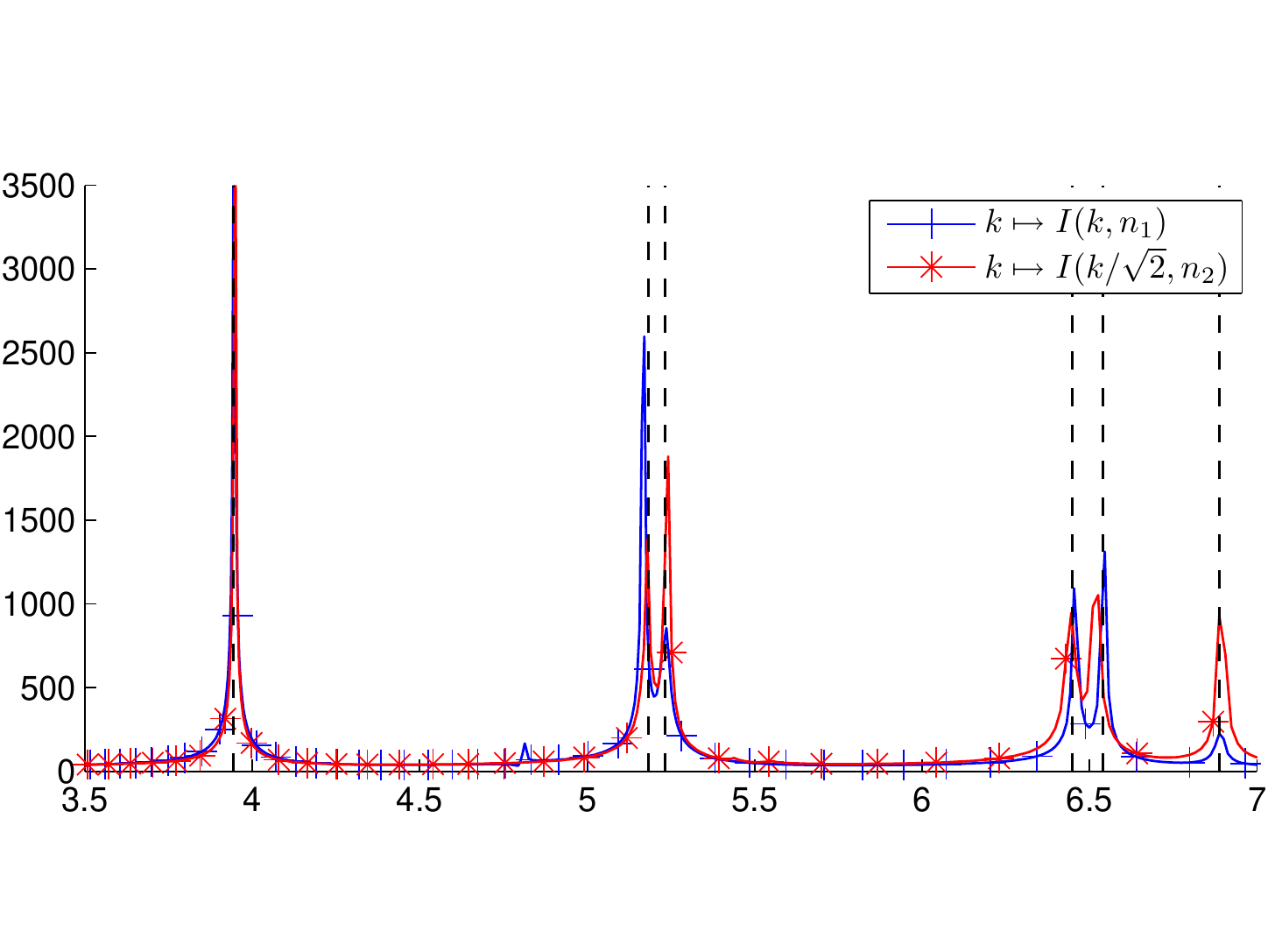}\\[-220pt]
\hspace{-5.1cm}\includegraphics[height=3.1cm]{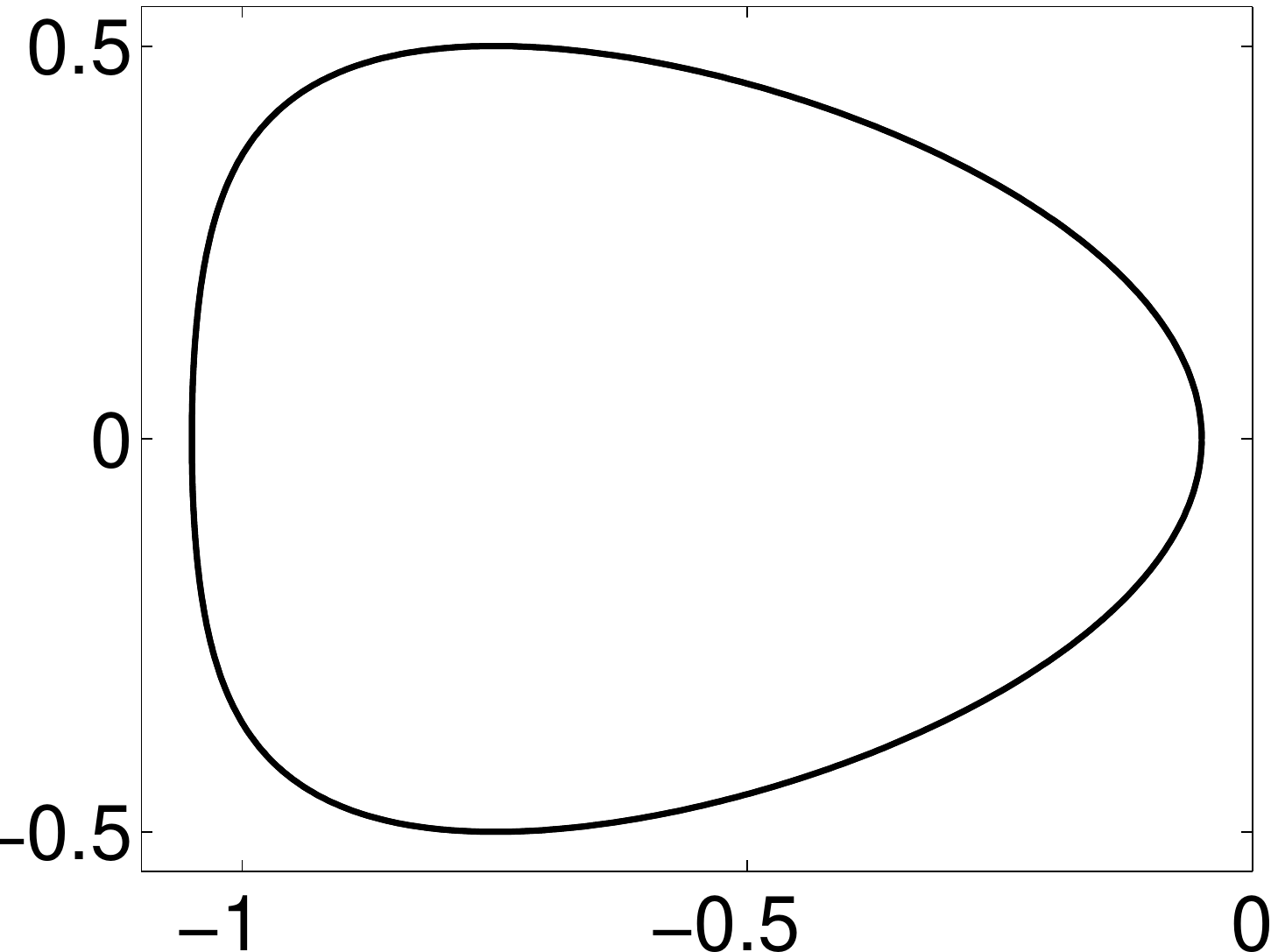}\\[130pt]
\end{tabular}
\caption{Curves $k\mapsto I(k,n_1)$ and  $k\mapsto I(k/\sqrt{2},n_2)$ for $k\in(3.5;7)$. The vertical dashed lines correspond to the eigenvalues of Problem \eqref{eigenpbSndOrder} computed with $n=n_1$ in $\Om$ (the quantities we want to retrieve in practice). The kite shaped domain represents the inclusion $\Om$.  \label{numerics}}
\end{figure}

We provide some preliminary numerical examples showing how the algorithm of determining constant $n$ using the ZIM background would work. The inclusion $\Om$ is a kite shape as depicted in Figure \ref{numerics}. First, we generate the farfield matrices for 200 incident directions $\boldsymbol{\theta}_{i}$ and 200 observation directions $\boldsymbol{\theta}_{s}$ by solving Problem (\ref{PbChampTotalFreeSpaceIntrod}). This is done for a range of wavenumbers $k$. Thus we obtain a discretization of the farfield operator $F$ defined in (\ref{defInitialFFoperator}) which in practice would be given by measurements. We consider two cases: $n=n_1=2.0$ in $\Om$ and $n=n_2=4.0$ in $\Om$. Then we choose $\Om_b=\Om$ and we build a discretization of the artificial farfield operator $F^{\mrm{art}}$ introduced in (\ref{DefArtFFop}). This allows us to compute the indicator function $I(k,n)$ defined in (\ref{indicatorGLSM}) with $\Om_0=\Om$ (for more details concerning this step, we refer the reader to \cite{l-thesis}). In Figure \ref{numerics}, we display the curve $k\mapsto I(k,n_1)$ for $k\in(3.5;7)$. On the other hand, using the \texttt{FreeFem++} software \cite{FreeFem}, we solve the eigenvalue Problem \eqref{eigenpbSndOrder} for $n=n_1$ in $\Om$. We observe that the square roots of the eigenvalues of \eqref{eigenpbSndOrder}, marked by dashed lines on Figure \ref{numerics}, correspond to the $k$ for which $k\mapsto I(k,n_1)$ has some peaks. This is coherent with the result of Theorem \ref{ThmCaractSpectrum}. Therefore, from the measurements, we can assess the $\lambda_p(n)$ in (\ref{totoratio}). Computing $\lambda_p(1)$ (solving \eqref{eigenpbSndOrder} with $n=1$), we can recover $n$ when $n$ is constant using Formula (\ref{totoratio}). In Figure \ref{numerics}, we also display the curve $k\mapsto I(k/\sqrt{2},n_2)$ for $k\in(3.5;7)$. We see that the peaks of this curve coincide with the ones of $k\mapsto I(k,n_1)$. 
This is in accordance with the theory which guarantees that $n_2/n_1=\lambda_p(n_1)/\lambda_p(n_2)$ for all $p\ge0$. Numerically $|n_2/n_1-\lambda_p(n_1)/\lambda_p(n_2)|$ is of order $10^{-2}$.

\end{document}